\documentclass[graybox]{svmult}

\usepackage{mathptmx}       
\usepackage{helvet}         
\usepackage{courier}        
\usepackage{type1cm}        
\usepackage{makeidx}         
\usepackage{graphicx}        
\usepackage{multicol}        
\usepackage[bottom]{footmisc}

\usepackage{mathtools}
\usepackage{amsfonts} 
\usepackage{amsmath}
\usepackage{bbold}
\usepackage[bb = boondox]{mathalpha}

\usepackage{stix}
\usepackage{amsmath,mleftright,mathtools}
\DeclareMathAlphabet\mathbfcal{OMS}{cmsy}{b}{n}

\usepackage{dsfont}
\usepackage[authoryear,round]{natbib}

\newcommand{\B}{\mathbb{B}}

\newcommand{\F}{\mathbb{F}}

\newcommand{\N}{\mathbb{N}}

\newcommand{\Q}{\mathbb{Q}}
\newcommand{\R}{\mathbb{R}}
\renewcommand{\S}{\mathbb{S}}

\newcommand{\X}{\mathbb{X}}

\newcommand{\CC}{\mathcal{C}}

\newcommand{\NN}{\mathcal{N}}

\newcommand{\PP}{\mathcal{P}}

\newcommand{\CCC}{\mathscr{C}}

\newcommand{\RRR}{\mathscr{R}}

\DeclareMathOperator{\Tr}{Tr}

\def\ppm{\scriptscriptstyle\pm} 
\def\F{F^{\ppm}}

\def\Q{Q^{\ppm}}

\def\FP{F^{\ppm}_{\rm P}}
\def\FPb{{\bf F}^{\ppm}_{\rm P}}
\def\QP{Q^{\ppm}_{\rm P}}
\def\QPb{{\bf Q}^{\ppm}_{\rm P}}
\def\R{{\mathbb R}}
\def\1{{\mathbb 1}}

\def\FPg{{\bf F}^{\scriptscriptstyle{\text{\rm g}}}_{\rm P} }
\def\QPg{{\bf Q}^{\scriptscriptstyle{\text{\rm g}}}_{\rm P} }

\def\g{
	\scriptscriptstyle{\text{\rm g}}
}

\makeindex             

\begin{document}

\title*{Multivariate Quantiles: Geometric and Measure-Transportation-Based Contours}
\titlerunning{ }
\author{Marc Hallin and Dimitri Konen}
\institute{Marc Hallin and Dimitri Konen \at Universit\' e libre de Bruxelles,  Belgium, and University of Warwick, United Kingdom
 \\\email{mhallin@ulb.be, dimitri.konen@warwick.ac.uk}
 }
%
%
\maketitle


\abstract{Quantiles are a fundamental concept in probability and theoretical statistics and a daily tool in their applications. While the univariate concept of quantiles is quite clear and well understood, its multivariate extension is more problematic. After half a century of continued efforts and many proposals, two concepts, essentially, are emerging: the so-called {\it (relabeled) geometric quantiles}, extending the characterization of univariate quantiles as minimizers of an L$_1$ loss function involving the {\it check functions}, and the more recent {\it center-outward quantiles} based on measure transportation ideas. These two concepts yield distinct families of {\it quantile regions} and {\it quantile contours}. Our objective here is to present a comparison of their main theoretical properties and a numerical investigation of their differences.}

\section{Introduction}
Quantiles are a fundamental concept  in probability and theoretical statistics and a familiar tool in their applications, ranging from descriptive statistics and data analysis to statistical inference. The univariate concept of a quantile function (the inverse of a distribution function) is well understood and well studied in dimension $d=1$. Despite half a century of efforts, however, the multivariate extension   (dimension $d\geq 2$) of this essential notion remains quite problematic---the obvious difficulty being the absence of a canonical ordering of the real space for dimension two and higher.  

Various proposals have been made in the literature, none of which is fully agreed upon. Two concepts, essentially, are emerging from these  many attempts: the so-called {\it (relabeled) geometric quantiles} (Chaudhuri 1996), extending the chara\-cte\-rization of univariate quantiles as minimizers of an L$_1$ loss function involving the {\it check functions}, and the more recent {\it center-outward quantiles} (Hallin et al.~2021) based on measure transportation ideas. 

The objective of this note is a brief comparison of some of the properties of the quantile regions and quantile contours associated with these two concepts.

\section{Quantiles: from univariate to multivariate} 

Throughout,  the inner product of two vectors $x,y\in\R^d$ is denoted as $x'y$; $\B_d$ and~$\S^{d-1}$ stand for the open unit ball and the $(d-1)$-dimensional unit sphere, respectively, in~$\R^d$.

\subsection{Definitions and basic properties}\label{11}

Let  $X$ denote a Lebesgue-absolutely continuous real-valued random variable with distribution $\rm P$ over $\R$ and distribution function 
\[
F_{\rm P}: x\mapsto F_{\rm P}(x)\coloneqq {\rm P}(X\leq x).
\]
To simplify the exposition,  assume that $X$ admits a nonvanishing density~$f\coloneqq {\rm dP/d\lambda}$ with respect to the Lebesgue measure $\lambda$ so that, in particular, $F_{\rm P}$ is continuous and  strictly monotone increasing, hence invertible. The quantile function of $\rm P$ (equivalently, the quantile function of $X\sim \rm P$) then is defined as the inverse $Q_{\rm P}\coloneqq F_{\rm P}^{-1}$ of the distribution function $F_{\rm P}$ and the quantile of order~$\alpha\in (0,1)$ of $\rm P$ (of $X\sim{\rm P}$) is the value~$Q_{\rm P}(\alpha)$ of $Q_{\rm P}$ at $\alpha$. 

Along with quantiles of given order $\alpha\in (0,1)$, it is often convenient to consider central quantile regions and contours. In dimension $d=1$,  
call {\it quantile region of order $\tau\in~\![0,1)$} of $\rm P$  (of $X\sim{\rm P}$) the closed interval
\[
\RRR_{\rm P}(\tau)
\coloneqq 
\Big[Q_{\rm P}\Big(\frac{1-\tau}{2}\Big),\, Q_{\rm P}\Big(\frac{1+\tau}{2}\Big)\Big]
\]
 (an interquantile interval) and  {\it quantile contour of order $\tau$ of $\rm P$}  (of $X\sim{\rm P}$) the boundary 
 \[
{\CCC_{\rm P}(\tau)}
\coloneqq  
\Big\{Q_{\rm P}\Big(\frac{1-\tau}{2}\Big),\, Q_{\rm P}\Big(\frac{1+\tau}{2}\Big)\Big\}
 \]
  of $\RRR_{\rm P}(\tau)$. 
  
  The quantile region of order 0 is the degenerate interval~ $\RRR_{\rm P}({0})=\{Q_{\rm P}(1/2)\}$ and has Lebesgue measure zero: call it the {\it median region}. Similarly, the quantile region of order 1/2 is the traditional  interquartile interval: $\RRR_{\rm P}({1/2})=[Q_{\rm P}(1/4), Q_{\rm P}(3/4)]$. Contrary to the quantile half\-lines~$(-\infty, Q_{\rm P}(\alpha)]$, which very much depend on the orientation of the real line (an orientation that is no longer meaningful in~${\mathbb R}^d$ for~$d\geq~\!2$), the concepts of quantile regions and contours of order $\tau$ are {\it ``center-outward''} concepts (the ``center'' being the median~$Q_{\rm P}(1/2)$) and are invariant under a change of the orientation of~$\R$ (an orthogonal transformation). The collection~$\left\{ \RRR_{\rm P}(\tau):\, \tau\in [0,1)
\right\}$ of quantile regions is stricly nested (in the sense of set inclusion) as $\tau$ ranges from zero to one, with~$\RRR_{\rm P}({0})=\bigcap_{\tau \in[0,1)}{\RRR_{\rm P}}(\tau)$. It is easily seen that ${\CCC_{\rm P}(\vert 2F_{\rm P}(x) -1\vert )}$ is the unique quantile contour  running through  $x\in\mathbb{R}$; actually,
\begin{equation}\label{eq1}
{\RRR_{\rm P}(\tau)} = \left\{ x : \vert 2F_{\rm P}(x) -1\vert \leq \tau
\right\}
 \ \text{ and }\  
{\CCC_{\rm P}(\tau)} = \left\{ x : \vert 2F_{\rm P}(x) -1\vert =\tau
\right\}
\quad \tau\in[0,1).
\end{equation}

An essential property of a quantile function is that the $\rm P$-probability content of a quantile region~$\RRR_{\rm P}(\tau)$ is $\tau$ irrespective of $\rm P$: namely, ${\rm P}(\RRR_{\rm P}(\tau)) = \tau$ for all $\rm P$. What would be the relevance, indeed, of an interquartile region ${\RRR_{\rm P}(1/2)}$ with probability content 0.5 under ${\rm P=P}_1$ (as should be) but  probability content 0.6 under ${\rm P=P}_2$, and  probability content 0.7 under ${\rm P=P}_3$? This property, in dimension $d=1$, is  equivalent to  the property of $F_{\rm P}(X)$ being uniform over $[0,1)$ for $X\sim {\rm P}$. 

When the center-outward quantile regions $\RRR_{\rm P}(\tau)$ and contours $\CCC_{\rm P}(\tau)$ are the main points of interest, it is natural to replace $F_{\rm P}$ and $Q_{\rm P}$ by  center-outward counterparts: call 
 $\FP\coloneqq 2F_{\rm P}-1$ and $\QP\coloneqq (\FP)^{-1}$ the {\it center-outward distribution} and {\it quantile functions} of $\rm P$ (or $X\sim{\rm P}$), respectively. Just like $F_{\rm P}$ and~$Q_{\rm P}$, the maps $\FP$ and~$\QP$ are monotone increasing and carry the same information about~$\rm P$, which they fully characterize. Instead of being uniformly distributed over~$[0,1)$,
 ~$\FP (X)$, when $X\sim{\rm P}$,  
  is uniform over $(-1, 1)$, which is the open unit ball $\B_1$ in~$\R$. In fact, denoting by ${\rm U}_1$ the uniform distribution over $\B_1$, we have $\F_{\rm P}(X)\sim U_1$ when~$X\sim{\rm P}$, hence $\QP (U)\sim{\rm P}$  if~$U\sim{\rm U}_1$. In the terminology of measure transportation, we say that the {\it transport}~$\FP$ is {\it pushing $\rm P$ forward to ${\rm U}_1$} and that its inverse~$\QP$ is {\it pushing ${\rm U}_1$ forward to ${\rm P}$}, which we conveniently denote as
 \[\FP\#{\rm P} = {\rm U}_1\quad\text{and }\quad \QP\#{\rm U}_1 = {\rm P}
 ,\]
 respectively. The quantile region $\RRR_{\rm P}(\tau)$ and quantile contour $\CCC_{\rm P}(\tau)$
 then take the simple forms (equivalent to \eqref{eq1})
 \begin{equation}\label{eq2}
 \RRR_{\rm P}(\tau)
 =
 \QP (\tau\B_1)
 \quad\text{and}\quad 
 \CCC_{\rm P}(\tau)
 =
 \QP (\tau\S^0)
 ,\quad \tau\in[0,1), 
 \end{equation}
respectively, where $\S^0\coloneqq \{-1, 1\}$ stands for the ($0$-dimensional) unit sphere in $\R$. 

\subsection{Characterization of  univariate quantiles as   minimizers of expected check functions}\label{12}

When $\rm P$ is a probability measure over $\R$ and $X\sim {\rm P}$, and given a level $\alpha\in [0,1)$, a quantile of order $\alpha$ of $X$, by definition, is any $x\in\R$  satisfying  
$$
{\rm P}(X<x)\leq \alpha \leq {\rm P}(X\leq x)
.
$$
It turns out that such $x$'s are exactly the minimizers (with respect to $\theta\in \R$) of the objective function 
\begin{equation}\label{eq:MinSpQuant}
\theta\mapsto M_\alpha^{\rm P}(\theta)
:=
{\rm E}_{\rm P}[\rho_\alpha(\theta-X)-\rho_\alpha(X)] 
\end{equation}
 where, denoting by $z^+$ and $z^-$ the positive and negative parts, respectively,  of $z\in\mathbb R$, 
$$
\rho_\alpha : z\mapsto
\rho_\alpha(z)
\coloneqq 
|z|+(2\alpha-1)z
=
2\big\{\alpha z^+ + (1-\alpha)z^-\big\}
,
\quad 
z\in\R,\ \alpha\in [0,1)
$$
stands for the so-called {\it check function}. Noting that $M_\alpha^{\rm P}(\theta)$ rewrites as
\begin{equation}\label{eq:GeometricFunction}
M_\alpha^{\rm P}(\theta)
=
{\rm E}_{\rm P}\big[|\theta-X|-|X|\big] - (2\alpha-1)\theta 
\end{equation}
for all $\theta\in\R$, it follows from the triangular inequality that  $\big\vert |X-\theta| - |X| \big\vert \leq |\theta|$, which is integrable.  The objective function~$M_\alpha^{\rm P}$, hence the correspondence between quantiles and the minimizing procedure in (\ref{eq:MinSpQuant}), thus remains well defined for an {\it arbitrary} probability measure ${\rm P}$, without any moment assumption. In this general setting, we  have,  even in the absence of a finite moment of order one,  
\[Q_{\rm P}(\alpha)={\rm argmin}_{\theta\in\R} M_\alpha^{\rm P}(\theta),\qquad \alpha\in [0,1).\]

 For $\alpha=0.5$, 
 $M^{\rm P}_\alpha(\theta)={\rm E}_{\rm P}[|\theta-X|-|X|]$; in case ${\rm E_P}(X)<\infty$, thus, we recover the familiar characterization of the median~$Q_{\rm P}(1/2)$ as an ${\rm L}_1$ location parameter, the minimizer of the expected absolute deviation ${\rm E}_{\rm P}[|\theta-X|]$.

 Recalling that $F_{\rm P}^{\ppm}\coloneqq 2F_{\rm P}-1$, we see that the natural labeling 
for  center-outward  quantiles is 
 provided by $\beta:=2\alpha-1\in(-1,1)$. Therefore, the center-outward version of the characterization of quantiles as minimizers is obtained by relabeling $2\alpha-1$ as $\beta$ in~\eqref{eq:MinSpQuant}. This leads to the characterization of center-outward quantiles $Q_{\rm P}^{\ppm} (\beta)$  as minimizers with respect to $\theta\in\R$ of the objective function 
 \begin{equation}\label{eq:GeometricFunctionRescale}
 \theta\mapsto O_{\beta}^{\rm P}(\theta)
 \coloneqq 
 {\rm E}_{\rm P}\big[|\theta-X|-|X|\big]-\beta \theta  ,\qquad \beta\in (-1,1).
 \end{equation}
 We then have (still, without any moment assumption)
 \[
 Q_{\rm P}^{\ppm}(\beta)={\rm argmin}_{\theta\in\R}O_\beta^{\rm P}(\theta), \qquad \beta\in (-1,1).
 \]
  The value $\beta=0$ yields the median $\QP(0)=Q_{\rm P}(1/2)$ of $\rm P$, while positive (resp. nega\-tive) values of $\beta$ 
   are associated to quantiles sitting to the right-hand (resp. left-hand)  side of the median.

\subsection{Characterization  of  univariate quantile functions as   monotone transports}\label{13}

Consider a probability measure ${\rm P}$ with  non-vanishing density $f$ with respect to the Lebesgue measure on $\R$. An alternative to the  characterization, in dimension $d=1$, of  $\QP$ (hence of $Q_{\rm P}$ and $\FP$)  as   minimizers of expected check functions is obtained by noting that~$\QP$ is monotone increasing and pushes ${\rm U}_1$ forward to   $\rm P$.  This fully characterizes~$\QP$ up to a set of values with Lebesgue measure zero.\footnote{As we shall see, the existence and (almost everywhere) uniqueness, in dimension $d=1$, of a monotone non-decreasing mapping (that is, the derivative, hence the gradient, of a convex function) pushing ${\rm U}_1$  forward to  $\rm P$ is a particular case of a famous and  more general  measure transportation result by \cite{McCann1995ExistenceAU}.}  The center-outward quantile function  $\QP$, thus, can be characterized as the almost everywhere unique monotone map from the open unit ball $\B_1=(-1,1)$ to $\R$ such that $\QP\#{\rm U}_1 = {\rm P}$.

\section{Multivariate quantile functions, regions, and contours}\label{2}

Extending to dimension $d>1$ the characterizations of quantile functions developed in Section~\ref{12}  
 leads to the concept of {\it geometric quantiles} first proposed  by \cite{Chaudhuri96} while   the measure-transportation-based characterization developed in Section~\ref{13} leads to the concept of multivariate {\it center-outward quantiles} proposed in \cite{Hallin2020DistributionAQ}.

\subsection{Geometric quantiles and  geometric quantile contours}\label{21}

Consider a probability measure ${\rm P}$ supported on $\R^d$, with $d\geq 2$. Extending to\linebreak a~$d$-dimensional framework the characterization of univariate quantiles as minimizers of expected check functions is obtained by replacing, in (\ref{eq:GeometricFunctionRescale}), the absolute  values with Euclidean norms and $\beta\in (-1,1)=\B_1$ with ${\bf v}\in \B_d$. Letting~${\bf v}\eqqcolon\tau {\bf u}\in\B_d$\linebreak with $\tau= \Vert {\bf v}\Vert \in (0,1)$ and ${\bf u}={\bf v}/\Vert {\bf v}\Vert \in\S^{d-1}$, define the {\it geometric quantile} of order~$\tau >0$ in direction ${\bf u}={\bf v}/\|{\bf v}\|$ of $\rm P$ (of ${\bf X}\sim {\rm P}$) as an arbitrary minimizer, over~${\bf z}\in\R^d$, of the objective function 
\begin{equation}\label{Eq3.5}
{\bf z}\mapsto 
O_{\tau,{\bf u}}^{\rm P}({\bf z})
:=
{\rm E}_{\rm P}\big[\|{\bf z}-{\bf X}\|-\|{\bf X}\|\big] - \tau {\bf u}'{\bf z}
\end{equation}
 where ${\bf X}\sim {\rm P}$; throughout, boldface is used to stress that the variables take values in~$\R^d$. The objective function $O_{\tau, {\bf u}}^{\rm P}$ being convex, it can be shown that the geo\-metric~quantiles of ${\rm P}$ are unique as soon as $\rm P$ is not supported on a single line of~$\R^d$---a situation 
 which essentially boils down to a univariate setting: 
 see Theorem~1 in \cite{PaiVir2020} for a   discussion.  
 For $\tau = 0$, define the {geometric median} as the minimizer of ${\rm E}_{\rm P}\big[\|{\bf z}-{\bf X}\|-\|{\bf X}\|\big]$ which, still in view of the triangular inequality, exists without any moment assumption and is better known as the {\it Fr\' echet median}.    
 
  When uniqueness holds ($\rm P$ not supported on a line), we denote by $\QPg(\tau {\bf u})$ the geometric quantile of order $\tau\in [0,1)$ in direction $\bf u$ of ${\rm P}$   (of ${\bf X}\sim{\rm P}$) and by $\QPg({\boldsymbol 0})$   the geometric or Fr\' echet median, which can be interpreted as the geometric quantile  of order $\tau=0$ in any direction $\bf u$.  Here again, boldface is used to stress that~$\QPb$ takes values in $\R^d$. 

Strongly related to the gradient of $O_{\tau,{\bf u}}^{\rm P}$, we define the {\it geometric distribution function} $ \FPg$ of ${\rm P}$ as the mapping 
\begin{equation}\label{eq:GeometricCdf}
{\bf z}\mapsto \FPg({\bf z})
\coloneqq 
{\rm E}_{\rm P}\Big[\frac{{\bf z}-{\bf X}}{\|{\bf z}-{\bf X}\|}\1[{\bf X}\neq {\bf z}]\Big]
,\qquad 
{\bf z}\in\R^d
,
\end{equation}
where $\1[A]$ stands for the indicator function of  $A$. While the definition of $\QPg$ is motivated by the characterization of univariate center-outward quantiles as minimizers of an objective function, let us stress that the geometric distribution function in (\ref{eq:GeometricCdf}) is a natural analogue of the univariate center-outward distribution as well. Indeed, recalling that, for $d=1$,  $F_{\rm P}^{\ppm}\coloneqq 2F_{\rm P}-1$  and observing that $2 \1[X\leq z]-1={\rm sign}(z-X)$, one obtains 
\begin{equation}\label{Eq4}
F_{\rm P}^{\ppm}(z)
=
{\rm E}_{\rm P}[{\rm sign}(z-X)]
,\quad 
z\in\R,
\end{equation}
to which $\FPg$ in (\ref{eq:GeometricCdf}) reduces when $d=1$. The direction  $({\bf z}-{\bf X})/\|{\bf z}-{\bf X}\|$, thus, can be interpreted as a multivariate sign for ${\bf z}-{\bf X}\neq {\boldsymbol 0}$. 

Provided that ${\rm P}$ is non-atomic and is not supported on a single line of $\R^d$, Theorem~6.2 in \cite{KonPai1} entails that $\FPg$ is a homeomorphism between $\R^d$ and the open unit ball $\B_d$ of $\R^d$, with inverse~$\QPg$. Consequently, one also can define $\QPg$ as $\left(\FPg\right)^{-1}$ with $\FPg$ defined in \eqref{eq:GeometricCdf}, or define $\FPg$ as $\left(\QPg\right)^{-1}$ with $\QPg$ defined as the minimizer of $O_{\tau,{\bf u}}^{\rm P}$ in~\eqref{Eq3.5}. Further results about the regularity of $\FPg$ and $\QPg$ are provided in \cite{Kon2022}. 

With this, one can define {\it geometric quantile regions}  and {\it contours}, which are multivariate extensions of the univariate center-outward quantile regions and contours, by letting 
\begin{equation}\label{eq3}
\pmb{\RRR}_{{\rm P}}^{\g}(\tau)
:=
\QPg(\tau \B_d)
\quad \textrm{and}\quad 
\pmb{\CCC}_{\rm P}^{\g}(\tau)
:=
\QPg(\tau\S^{d-1})
,\quad 
\tau\in [0,1)
.
\end{equation}
For $\tau = 0$,   $\pmb{\RRR}_{\rm P}^{\g}(0)$ and $\pmb{\CCC}_{\rm P}^{\g}(0)$ reduce to the singlepoint set consisting of the {\it geometric} or  {\it Fréchet median}; as in the univariate case, the geometric median is an ${\rm L}_1$ location parameter that minimizes, when finite first-order moments exist, the expected absolute deviation ({\it absolute} deviations, here, are to be understood as {\it Euclidean norms} of deviations).

\subsection{Measure-transportation-based quantiles and  quantile contours}\label{22}

Denote by ${\rm U}_d$ the spherical uniform over the open unit ball $\B_d$, that is, the product of a uniform distribution over the unit sphere $\S^{d-1}$ and a (univariate) uniform over the distances from the origin; in other words, ${\rm U}_d$ is the probability measure with density~$f_{{\rm U}_d}({\bf x})$ proportional to $1/\|{\bf x}\|^{d-1}\1[{\boldsymbol 0}<\|{\bf x}\|\leq 1]$ for ${\boldsymbol 0}<\|{\bf x}\|$ and  $f_{{\rm U}_d}({\boldsymbol 0})=\infty$.

A celebrated result by \cite{McCann1995ExistenceAU} implies that, for any distribution~$\rm P$ over~${\R} ^d$, there exists an almost everywhere unique gradient of a convex function pushing $\rm P$  forward to the spherical uniform ${\rm U}_d$ over the (open) unit ball $\B_d$ of $\R^d$. In dimen\-sion~$d=1$,  a gradient of convex function is a monotone non-decreasing function and the unique monotone non-decreasing function pushing~$\rm P$ forward to ${\rm U}_1$ is the center-outward distribution function $\FP\coloneqq 2F_{\rm P}-1$. For  arbitrary dimension $d\geq~\!1$, this leads~\cite{Hallin2020DistributionAQ} to define the {\it center-outward distribution function} of~$\rm P$ (or~${\bf X}\sim{\rm P}$) as the almost everywhere unique gradient of a convex function~$\FPb$ pushing~$\rm P$  to the spherical uniform ${\rm U}_d$.

Assuming  that $\rm P$ belongs to the class $\mathcal{P}_d^+$ of  Lebesgue-absolutely continuous distributions with density bounded away from zero and infinity on any compact\footnote{Formally, for any compact subset $K\subset \R^d$, there exist constants $0<c^-_K\leq c_K^+<\infty$ such \linebreak that~$c_K^-\leq {\rm dP}/{\rm d\lambda}({\bf x})\leq c_K^+$ for all ${\bf x}\in K$. This assumption is maintained here for clarity of exposition, but it can be relaxed: see \cite{BarrioSanzHal} and  \cite{delGonz2023}.} subset of~$\R^d$,  \cite{FigalliCenter} showed that $\FPb$ is a homeo\-morphism between~$\mathbb{R}^d\!\setminus\!(\FPb)^{-1}({\bf 0})$ and the punctured unit ball~$\B_d\!\setminus\!\{{\bf 0}\}$. 
 The restriction of $\FPb$ to~$\mathbb{R}^d\!\setminus\!(\FPb)^{-1}({\bf 0})$ is thus  continuously invertible. We call its (continuous) inverse $\QPb\coloneqq (\FPb)^{-1}$ the {\it center-outward quantile function} of $\rm P$ (or ${\bf X}\sim{\rm P}$). This yields  {\it center-outward quantile regions}  and {\it contours}, of order $\tau$ which are multivariate extensions of the univariate center-outward quantile regions and contours defined in~\eqref{eq2}, of the form %
%
 \begin{equation}\label{eq4}
\pmb{\RRR}_{\rm P}^{\ppm}(\tau)
\coloneqq 
\QPb (\tau \B_d)
\quad\text{and}\quad 
\pmb{\CCC}_{\rm P}^{\ppm}(\tau)
\coloneqq 
\QPb (\tau \S^{d-1})
,\qquad  \tau\in(0,1),
 \end{equation}
respectively. As for $\tau = 0$, define~$\pmb{\RRR}_{\rm P}^{\ppm}(0)\coloneqq \bigcap_{0<\tau <1}\pmb{\RRR}_{\rm P}^{\ppm}(\tau)$ as the {\it center-outward median region}, which is shown to be convex and compact, with Hausdorff dimension at most $d-1$ \citep{FigalliCenter}. Note the parallel between \eqref{eq4} and \eqref{eq3}. 

\section{Main properties}

\subsection{Geometric quantile regions and contours} 

Consider a probability measure ${\rm P}$ which is non-atomic and is not supported on a single line of $\R^d$; in particular, $\FPg$ and $\QPg$ are homeomorphisms with $\QPg = (\FPg)^{-1}$. Then, geometric quantile regions and contours have the following topological properties:

\begin{enumerate}
	\item[(i)] the regions $\pmb{\RRR}_{\rm P}^{\g}(\tau)$ are compact and arc-connected for all $\tau\in[0,1)$;
	
	\item[(ii)] the regions $\pmb{\RRR}_{\rm P}^{\g}(\tau)$ are stricly nested as $\tau$ increases: for any~$0\leq\tau_1 < \tau_2<1$, we have $\pmb{\RRR}_{\rm P}^{\g}(\tau_1) \subsetneq \pmb{\RRR}_{\rm P}^{\g}(\tau_2)$; 
	
	
	\item[(iii)] if ${\rm P}$ has a bounded density $f$, then the contour $\pmb{\CCC}_{\rm P}^{\g}(\tau)$ is a $(d-1)$-dimensional manifold of class $\CC^{d-1}$ for all $\tau \in (0,1)$; if, in addition, $f$ is of class $\CC^{k,\alpha}$ for some~$k\in\N$ and $\alpha\in (0,1)$, then $\pmb{\CCC}_{\rm P}^{\g}(\tau)$ is of class $\CC^{d+k}$ for all $\tau\in (0,1)$;

	\item[(iv)] for $d=1$, $\FPg$, $\QPg$, $\pmb{\RRR}_{\rm P}^{\g}$, and $\pmb{\CCC}_{\rm P}^{\g}$ reduce to $\FP$, $\QP$, ${\RRR}_{\rm P}^{\ppm}$, and ${\CCC}_{\rm P}^{\ppm}$. 
\end{enumerate}

While properties (i)-(ii) are straightforward consequences of the continuity of~$\QPg$, property~(iii) requires more advanced mathematical tools. That property  was established in \cite{Kon2022}, where it is shown that geometric distribution functions sa\-tisfy a linear partial differential equation, which was the key to te derivation of the regularity properties of geometric contours. Among other things, this PDE implies that geometric distribution functions fully characterize {\it arbitrary} probability measures (no moment or non-atomicity assumptions required), namely, 
\begin{enumerate}
	\item[(v)] if ${\rm P_1}$ and ${\rm P_2}$ are {\it arbitrary} probability measures such that ${\bf F}_{\rm P_1}^{\g} = {\bf F}_{\rm P_2}^{\g}$, then ${\rm P}_1={\rm P}_2$.
\end{enumerate} 
This characterization result has been long known, though: the first proof was provided in Theorem 2.5 of \cite{Koltchinski97} and was refined recently  in \cite{Kon2022}.\\

Geometric quantile and distribution functions further enjoy the following natural equivariance properties; see \cite{Girard17} and \cite{KonPai1}. 

%
%
%
\begin{enumerate}
	\item[(vi)] Let $\boldsymbol{\theta}\in\R^d$ and denote by ${\bf O}$ an arbitrary $d\times d$ orthogonal matrix. For ${\bf X}\sim {\rm P}$, write ${\bf F}_{ {\bf X}}^{\g} , {\bf F}_{{\bf X}}^{\g}  $, etc. instead of   ${\bf F}_{\rm P}^{\g} , {\bf F}_{\rm P}^{\g}  $, etc.  We have 
	\begin{align*}
		&{\bf F}_{\boldsymbol{\theta}+ {\bf O}{\bf X}}^{\g} (\boldsymbol{\theta} + {\bf O}{\bf x})
		= 
		{\bf O}{\bf F}_{\bf X}^{\g}({\bf x})
		\quad 
		&{\bf x}\in\mathbb{R}^d,
		&\\
		&{\bf Q}_{\boldsymbol{\theta}+ {\bf O}{\bf X}}^{\g} ({\bf O}{\bf u})
		= 
		\boldsymbol{\theta} + {\bf O}{\bf Q}_{\bf X}^{\g}({\bf u})
		\quad 
		&{\bf u}\in\B_d,&\\
		&\pmb{\RRR}_{\boldsymbol{\theta}+ {\bf O}{\bf X}}^{\g}(\tau)
		= 
		\boldsymbol{\theta}+ {\bf O}\pmb{\RRR}_{{\bf X}}^{\g}(\tau)
		\quad 
		&\tau\in[0,1), &\quad \text{ and}
		\\
		&\pmb{\CCC}_{\boldsymbol{\theta}+ {\bf O}{\bf X}}^{\g}(\tau)
		= 
		\boldsymbol{\theta}+ {\bf O}\pmb{\CCC}_{\bf X}^{\g}(\tau)
		\quad 
		&\tau\in[0,1).
	\end{align*}
\end{enumerate}
In particular, if $\rm P$ enjoys symmetry properties (rotational, central, or with respect to some hyperplanes), $\FPg$, $\QPg$, $\pmb{\RRR}_{\rm P}^{\g}(\tau)$, and $\pmb{\CCC}_{\rm P}^{\g}(\tau)$ enjoy {\it the same} symmetry properties as $\rm P$ for all $\tau\in~\![0,1)$.\\


Recall that an essential property of distribution functions in dimension $d=1$ is their non-decreasing monotonicity. In higher dimensions, the concept of monotonicity still makes sense; we distinguish two notions of monotonicity, which both reduce to the usual concept of (non-decreasing) monotonicity when $d=1$: {\it monotonicity} and {\it cyclical monotonicity}. We say that a map ${\bf G}:\R^d\to\R^d$ is {\it monotone} if 
$$
({\bf G}({\bf y})-{\bf G}({\bf x}))'({\bf y}-{\bf x})\geq 0
$$
for all ${\bf x}, {\bf y}\in\R^d$. We say that $G$ is {\it cyclically monotone} if for any $m\geq 1$ we have
$$
\sum_{k=1}^m \big({\bf G}({\bf x}_{k+1})-{\bf G}({\bf x}_k)\big)'{\bf x}_{k+1}
\geq 
0 
$$
for every cycle $x_1, x_2,\ldots, x_{m+1}:=x_1$. For $m=2$, cyclical monotonicity yields 
$$
\big({\bf G}({\bf x}_2)-{\bf G}({\bf x}_1)\big)'{\bf x}_2 
+
\big({\bf G}({\bf x}_1)-{\bf G}({\bf x}_2)\big)'{\bf x}_1
\geq 
0 
,
$$
and cyclical monotonicity thus implies monotonicity. 
Theorem 24.8 in \cite{Rocka70} essentially entails that ${\bf G}$ is cyclically monotone precisely when it writes as the gradient~${\bf G}=\nabla \psi$ of some convex (potential)  function~$\psi:\R^d\to\R$. 

The following then holds.\begin{enumerate}
	\item[(vii)] The maps $\FPg$ and $\QPg$ are cyclically monotone, hence also monotone, over $\R^d$.
\end{enumerate}
Indeed, similarly to the center-outward distribution functions based on measure transportation introduced in Section 2.1, geometric distribution functions are gradients of convex functions. This is clear from (\ref{eq:GeometricCdf}), which immediately yields $\FPg = \nabla \phi_{\rm P}^{\g}$, where 
$${\bf z}\mapsto 
\phi_{\rm P}^{\g}({\bf z})
:=
{\rm E}\big[ \|{\bf z}-{\bf X}\| - \|{\bf X}\|\big]
,\quad 
{\bf z}\in\R^d
$$
is a convex function and ${\bf X}\sim {\rm P}$. It follows that~$\FPg$ is cyclically monotone, which implies that its inverse $\QPg$ is cyclically monotone too (hence the gradient  of a convex function).

However, the geometric distribution function ${\bf F}_{\rm P}^{\g}$,  in general, does not push ${\rm P}$ forward to the spherical uniform ${\rm U}_d$ over $\B_d$. Would this be case, $\FPg$ and $\FPb$ both being gradients of  convex functions, the a.s.\ uniqueness result of McCann~(1955)  would imply $\FPg = \FPb$ a.s.\ for all ${\rm P}\in \PP_d^+$ (see Section 2.2), which  cannot hold in general. Indeed, $\FPg$ can be defined as the solution of the previously mentioned {\it linear} PDE   (see \cite{Kon2022}), while $\FPb$ can be defined as the solution of a  PDE of the Monge-Amp\`ere type which  is {\it nonlinear} (see \cite{Hallin2020DistributionAQ}). Therefore, the solu\-tions~$ \FPg$ and $ \FPb$, for general $\rm P$,  do not coincide and $\FPg$,  in general,   does not push~${\rm P}$ forward to~${\rm U}_d$.

Nevertheless, recalling that $\FPg$ is a homeomorphism between $\R^d$ and $\B_d$, we see that each point ${\bf x}\in\R^d$  belongs to exactly one geometric contour, namely ${\bf x}\in\pmb{\CCC}_{\rm P}^{\g}(\tau)$ where $\tau =\Vert \FPg({\bf x})\Vert$. Therefore, we can relabel quantile regions so that, in the new labeling, the region associated with an order $\tau\in [0,1)$ has ${\rm P}$-probability content $\tau$. For this purpose, it is enough to observe that the map 
$$
\tau\mapsto \kappa_{\rm P}^{\g}(\tau)
:=
{\rm P}[\pmb{\RRR}_{\rm P}^{\g}(\tau)]
,\quad 
\tau\in [0,1),
$$ from $[0,1)\to [0,1)$ 
is a continous and monotone increasing bijection. In fact, $\kappa_{\rm P}^{\g}$ is the univariate distribution function of the random variable $\|\FPg({\bf X})\|$ with ${\bf X}\sim~\!{\rm P}$.\linebreak  Consequently, we can define the {\it relabeled geometric distribution function} of $\rm P$ (of~${\bf X}\sim {\rm P}$) as 
\begin{equation}\label{eq:RelabCdf}
{\bf x}\mapsto \tilde{\bf F}_{\rm P}^{\g}({\bf x})
:=
{\bf B_{\rm P}^{\g}}\big(\FPg({\bf x})\big)
,
\end{equation}
where 
$$
{\bf B_{\rm P}^{\g}}({\bf x})
:= 
\begin{cases} 
	\kappa_{\rm P}^{\g}(\|{\bf x}\|) \frac{{\bf x}}{\|{\bf x}\|} & \textrm{ if } {\bf x} \neq 0\\
	0 & \textrm{ if } {\bf x} = 0
\end{cases} 
$$
is a homeomorphism.
 Consequently, the relabeled geometric distribution function~$\tilde{\bf F}_{\rm P}^{\g}$ is also a homeomorphism, hence continuously invertible:   define the {\it relabeled geometric quantile function} as its inverse $\tilde{\bf Q}_{\rm P}^{\g}:=(\tilde{\bf F}_{\rm P}^{\g})^{-1}$. The collections of regions and contours associated with $\tilde{\bf F}_{\rm P}^{\g}$ are the same as those   associated with $\FPg$---only the labels have changed. 
 Call $
 	\tilde{\pmb{\RRR}}_{\rm P}^{\g}(\tau)
 	:=
 	\tilde{\bf Q}_{\rm P}^{\g}(\tau \B_d)
 	$ the {\it relabeled geometric quantile region of order $\tau$}.  Then,  $\tilde{\pmb{\RRR}}_{\rm P}^{\g}(\tau)$ satisfies the essential property of a quantile region of order~$\tau$:
 \begin{enumerate}
 	\item[(viii)] the relabeled geometric quantile region 
 	$
 	\tilde{\pmb{\RRR}}_{\rm P}^{\g}(\tau)
 	$
 	has ${\rm P}$-probability content $\tau$ for all~$\tau\in~\![0,1)$, irrespective of $\rm P$.
 \end{enumerate}  
We similarly define the {\it relabeled geometric quantile contours}  
 $$
 \tilde{\pmb{\CCC}}_{\rm P}^{\g}(\tau)
 :=
 \tilde{\bf Q}_{\rm P}^{\g}(\tau \S^{d-1})
 ,\quad \tau \in [0,1)
 .
 $$
 
 One may wonder whether $\tilde{\bf F}_{\rm P}^{\g}$ now coincides with $\FPb$ (hence, $\tilde{\bf Q}_{\rm P}^{\g}$  with $\QPb$). First, let us stress that, even though $\|\tilde{\bf F}_{\rm P}^{\g}({\bf X})\|$, just like $\|\FPb({\bf X})\|$,  is uniformly distributed over $[0,1)$ when ${\bf X}\sim {\rm P}$, the behavior of~$\FPg({\bf X})/\|\FPg({\bf X})\|$, to the best of our knowledge, remains unknown, and has not been investigated so far. Second, observe that ${\bf B}_{\rm P}^{\g}$ writes as the gradient~$\nabla h_{\rm P}^{\g}$ of the convex function 
 $${\bf x}\mapsto 
 h_{\rm P}^{\g}({\bf x})
 :=
 \int_0^{\|{\bf x}\|} \kappa_{\rm P}^{\g}(s)\, ds 
 ,\quad 
 {\bf x}\in\R^d 
 .
 $$
 It follows that $\tilde{\bf F}_{\rm P}^{\g} = \nabla h_{\rm P}^{\g} \circ \nabla \phi_{\rm P}^{\g}$. However, it is well known that the composition of two gradients of convex functions, in general, is not the gradient of a convex function. For these reasons, which do not constitute rigorous proofs, we believe that there is no obvious reason why the equality $\tilde{\bf F}_{\rm P}^{\g}=\FPb$ should hold for general $\rm P$. The simulations in Section~\ref{4}) below---which do not constitute mathematical proofs---actually suggest that~$ \tilde{\pmb{\CCC}}_{\rm P}^{\g}(\tau)$   and~$\pmb{\CCC}_{\rm P}^{\ppm}(\tau)$, for some distributions~$\rm P$, very strongly differ. 
 

\subsection{Center-outward quantile regions and contours}

Let ${\rm P}\in\mathcal{P}_d^+$. It is shown in \cite{Hallin2020DistributionAQ} that 

\begin{enumerate}
\item[(i)]the regions  $\pmb{\RRR}_{\rm P}^{\ppm}(\tau)$ are compact and arc-connected for all $\tau\in(0,1)$\vspace{1mm};

\item[(ii)]the regions $\pmb{\RRR}_{\rm P}^{\ppm}(\tau)$ are stricly nested as $\tau$ increases:  $\pmb{\RRR}_{\rm P}^{\ppm}(\tau_1)\subsetneq\pmb{\RRR}_{\rm P}^{\ppm}(\tau_2)$ for \linebreak any~$0<\tau_1 < \tau_2<1$; 

\item[(iii)] the median region $\pmb{\RRR}_{\rm P}^{\ppm}(0)$  is compact, convex, and has Hausdorff dimension at most $(d-1)$;

\item[(iv)] if, for some $k\in\N$ and $\alpha\in (0,1)$, the density of ${\rm P}$ is locally of class $\CC^{k,\alpha}$, then the contour $\pmb{\CCC}_{\rm P}^{\ppm}(\tau)$ is a $(d-1)$-dimensional manifold of class $\CC^{k+1,\alpha}$ for all $\tau \in (0,1)$;\vspace{1mm}

\item[(v)] for $d=1$, $\FPb$, $\QPb$, $\pmb{\RRR}_{\rm P}^{\ppm}$, and $\pmb{\CCC}_{\rm P}^{\ppm}$ reduce to $\FP$, $\QP$, ${\RRR}_{\rm P}^{\ppm}$, and ${\CCC}_{\rm P}^{\ppm}$. \end{enumerate}

Note that (iv) was first established in \cite{BarrioSanzHal} when the support of~${\rm P}$ is convex, then refined in \cite{delGonz2023}. By construction, we also have the following:

\begin{enumerate}
\item[(vi)]the $\rm P$-probability content of $\pmb{\RRR}_{\rm P}^{\ppm}(\tau)$ equals the ${\rm U}_d$-probability content of $\tau\B_d$, which is $\tau$:  $\pmb{\RRR}_{\rm P}^{\ppm}(\tau)$ and $\pmb{\CCC}_{\rm P}^{\ppm}(\tau)$, thus, qualify as the center-outward quantile region and  quantile contour {\it of order} $\tau$ of $\rm P$;

\item[(vii)]if ${\bf X}\sim{\rm P}$, then $\FPb ({\bf X})\sim{\rm U}_d$; if ${\bf U}\sim{\rm U}_d$, then $\QPb ({\bf U})\sim{\rm P}$: both $\FPb$ and $\QPb$, thus, fully characterize $\rm P$.\footnote{Note, however, that the collections of center-outward  regions or   center-outward contours alone (that is, without the mappings $\FPb$ or $\QPb$ between $\R^d$ and $\B_d$) {\it do not} characterize $\rm P$.}

\item[(viii)] the maps $\FPb$ and $\QPb$ are cyclically monotone, hence monotone, over $\R^d$.
\end{enumerate}

Finally, $\FPb$ and $\QPb$ enjoy the following invariance/equivariance properties (see Hallin et al.~(2022)).

\begin{enumerate}
	\item[(ix)] Let $\boldsymbol{\theta}\in\R^d$ and  denote by ${\bf O}$ an arbitrary $d\times d$ orthogonal matrix. Letting ${\bf X}\sim {\rm P}$, write  ${\bf F}_{{\bf X}}^{\ppm}$, ${\bf Q}_{{\bf X}}^{\ppm}$, etc. for ${\bf F}_{{\rm P}}^{\ppm}$, ${\bf Q}_{\rm P}^{\ppm}$, etc. Then, 
	\begin{align*}
		&{\bf F}_{\boldsymbol{\theta}+ {\bf O}{\bf X}}^{\ppm} (\boldsymbol{\theta} + {\bf O}{\bf x})
		= 
		{\bf O}{\bf F}_{\bf X}^{\ppm}({\bf x})
		\quad 
		&{\bf x}\in\mathbb{R}^d,
		&\\
		&{\bf Q}_{\boldsymbol{\theta}+ {\bf O}{\bf X}}^{\ppm} ({\bf O}{\bf u})
		= 
		\boldsymbol{\theta} + {\bf O}{\bf Q}_{\bf X}^{\ppm}({\bf u})
		\quad 
		&{\bf u}\in\B_d,&\\
		&\pmb{\RRR}_{\boldsymbol{\theta}+ {\bf O}{\bf X}}^{\ppm}(\tau)
		= 
		\boldsymbol{\theta}+ {\bf O}\pmb{\RRR}_{{\bf X}}^{\ppm}(\tau)
		\quad 
		&\tau\in[0,1), &\quad \text{ and}
		\\
		&\pmb{\CCC}_{\boldsymbol{\theta}+ {\bf O}{\bf X}}^{\ppm}(\tau)
		= 
		\boldsymbol{\theta}+ {\bf O}\pmb{\CCC}_{\bf X}^{\ppm}(\tau)
		\quad 
		&\tau\in[0,1).
	\end{align*}
\end{enumerate}
In particular, if $\rm P$ enjoys symmetry properties (rotational, central, or with respect to some hyperplanes), $\FPb$, $\QPb$, $\pmb{\RRR}_{\rm P}^{\ppm}(\tau)$, and $\pmb{\CCC}_{\rm P}^{\ppm}(\tau)$ enjoy the same symmetry properties as $\rm P$ for all $\tau\in~\![0,1)$. This is not the case for the distribution and quantile functions resulting from a transport to the Lebesgue uniform  over the unit cube $[0,1]^d$ which are favored by some authors.\\

{These properties, which  for dimension $d=1$ reduce to those of $\F$ and $\Q$ and constitute the essence of the concepts of distribution and quantile functions, fully justify the terminology {\it center-outward distribution and quantile functions} for $\FPb$ and $\QPb$.}  {The {\it only} difference, here, between the properties of~$\FPb$ and~$\QPb$ and those of~$\FPg$ and~$\QPg$ is~(vi) above,} and this difference disappears if the relabeled functions~$\tilde{\bf F}_{\rm P}^{\g}$ and~$\tilde{\bf Q}_{\rm P}^{\g}$ are substituted for ${\bf F}_{\rm P}^{\g}$ and~${\bf Q}_{\rm P}^{\g}$, respectively.

In view of this,  the relabeled geometric and the center-outward concepts, which, it seems, do not coincide, apparently both qualify as equally satisfactory multivariate concepts of distribution and quantile functions. Is there a way to break that tie and select one of them rather than the other? Since analytical expressions are impossible, one has to recur to numerical approximations---in statistical terms, simulations---in order to perform comparisons. This is what we are doing in the next section, where comparisons are based on  simulations of the quantile contours $\pmb{\CCC}_{\rm P}^{\ppm}(\tau)$ and $\pmb{\CCC}_{\rm P}^{g}(\tau)$.

\section{Geometric 
 and center-outward quantile contours: an empirical comparative study}\label{4}

\subsection{Numerical computation of contours}\label{41}

Except for some very special distributions, such as the spherical ones, geometric and center-outward quantiles in dimension $d\geq 2$, unfortunately,  are sharing the inconvenient absence of  analytical expressions. Numerical solutions are possible, though, based on  simulations with high numbers of replications. For all $N\in\mathbb N$, denote by~${\rm P}_N$ the empirical distribution of a random sample of size $N$ drawn from ${\rm P}$. The accuracy of these numerical solutions relies on  consistency results for the empirical counter\-parts~${\bf F}_{{\rm P}_N}^{\g}$ and ${\bf Q}_{{\rm P}_N}^{\g}$ in the geometric case (obtained by substituting ${\rm P}_N$ for ${\rm P}$ in (\ref{eq:GeometricCdf})), for ${\bf F}_{{\rm P}_N}^{\ppm}$ and ${\bf Q}_{{\rm P}_N}^{\ppm}$ in the measure-transportation setting (see \cite{Hallin2020DistributionAQ} for precise definitions and details). 

In the geometric case, \cite{Mottonen97} established a Glivenko-Cantelli result of the form 
\begin{equation}\label{GCg}
\sup_{{\bf z}\in\R^d} 
\| {\bf F}_{{\rm P}_N}^{\g}({\bf z}) - \FPg({\bf z})\| 
\to 
0 
\end{equation}
with ${\rm P}$-probability one as $N\to\infty$, provided that ${\rm P}$ admits a bounded density. In fact, one can show that the bounded density assumption is superfluous: similarly to distribution functions in dimension $d=1$, the  Glivenko-Cantelli result \eqref{GCg} is valid for an {\it arbitrary} probability measure. The same convergence holds for the  sequence of relabeled counterparts $\tilde{{\bf F}}_{\rm P}^{\g}$ of ${{\bf F}}_{\rm P}^{\g}$ and $\tilde{{\bf F}}_{{\rm P}_N}^{\g}$ of ${{\bf F}}_{{\rm P}_N}^{\g}$; see \cite{Kon2022}.

In the  measure transportation approach, \cite{Hallin2020DistributionAQ} proved a similar Glivenko-Cantelli result for an interpolated version $\overline{{\bf F}}_{{\rm P}_N}^{\,\ppm}$  of their empirical center-outward distribution functions ${\bf F}_{{\rm P}_N}^{\ppm}$:
$$
\sup_{{\bf z}\in\R^d} 
\| \overline{{\bf F}}_{{\rm P}_N}^{\,\ppm}({\bf z}) - \FPb({\bf z})\| 
\to 
0 
$$
with ${\rm P}$-probability one as $N\to\infty$, when ${\rm P}\in\PP_d^+$ (see Section 2.2 for the definition of $\PP_d^+$).

\subsection{Numerical comparisons of quantile contours}\label{42}

Figure 1 shows, for orders $\tau\in \{.25, .50, .75\}$,  plots of the empirical contours $\tilde{\pmb{\CCC}}_{{\rm P}_N}^{\g}(\tau)$\linebreak  and~$\pmb{\CCC}_{{\rm P}_N}^{\ppm}(\tau)$ based on  random samples $\X _N=\{{\bf X}_1,\ldots, {\bf X}_N\}$ of size $N=2,400$ from four different probability distributions on $\R^2$. Given the regular grid 
$$
\Big\{ 
	{\bf u}_k = \big(\cos(\theta_k), \sin(\theta_k)\big) : \theta_k = \frac{2k\pi}{70} \textrm{ and } 1\leq k\leq K
\Big\} 
$$
of directions on $\S^1$, with $K=70$, each contour $\tilde{\pmb{\CCC}}_{{\rm P}_N}^{\g}(\tau)$ is estimated by 
$$
\Big\{ \tilde{{\bf Q}}_{{\rm P}_N}^{\g}(\tau {\bf u}_k) : 1\leq k \leq K=70 \Big\}
,
$$
where $\tilde{{\bf Q}}_{{\rm P}_N}^{\g}$ stands for the relabeled geometric quantile map defined in Section 3.1. In other words, to each level $\tau$ we first associate the level $\tau'$ such that 
$$
\tilde{\pmb{\CCC}}_{{\rm P}_N}^{\g}(\tau) = \pmb{\CCC}_{{\rm P}_N}^{\g}(\tau')
.
$$
Then, the geometric contour $\pmb{\CCC}_{{\rm P}_N}^{\g}(\tau')$ is obtained as the collection $\{{\bf z}_1,\ldots, {\bf z}_{K}\}$ such that ${\bf z_k}$ minimizes the empirical objective function $O_{\tau', {\bf u_k}}^{{\rm P}_N}$. 

Each measure-transportation-based contour is estimated by
$$
\Big\{ 
	{\bf Q}_{{\rm P}_N}^{\ppm}(\tau {\bf u}_k) : 1\leq k\leq K
\Big\}
.
$$
To construct the map ${\bf Q}_{{\rm P}_N}^{\ppm}$, we first fix a regular $n_R\times n_S$ grid ${\mathfrak G}_N$ over the punctured unit ball $(0,1)\times \S^1$ 
 such that $n_R  n_S = N$;\footnote{In the case of a genuine sample, $N$  is factorized into $N=n_Rn_S + n_0$ with $n_0< \min (n_R, n_S)$. We refer to \cite{HMord23} for details on the choice of $n_R$ and $n_S$. Here, however, we are dealing with simulations, and can choose $N$ such that $n_0=0$.} our simulation relies on~$N=2,400$ replications, with $n_R = 40$ and $n_S = 60$. Then, we seek an optimal coupling between~$\X _N$ and ${\mathfrak G}_N$, i.e., we associate to each ${\bf X}_k\in \X _N$ a unique ${\bf F}_{{\rm P}_N}^{\ppm}({\bf X}_k)\in \mathfrak G$ in such a way that
$$
\sum_{k=1}^N \|{\bf X}_k - {\bf F}_{{\rm P}_N}^{\ppm}({\bf X}_k)\|^2 
=
\min_{\pi\in\Pi} \sum_{k=1}^N \|{\bf X}_{\pi(k)} - {\bf F}_{{\rm P}_N}^{\ppm}({\bf X}_k)\|^2 
,
$$
where $\pi$ ranges over the set $\Pi$ of permutations of $\{1,\ldots, N\}$. A continuous map $\overline{\bf F}_{{\rm P}_N}^{\,\ppm}$ then is obtained from the values ${\bf F}_{{\rm P}_N}^{\,\ppm} ({\bf X}_k)$ of ${\bf F}_{{\rm P}_N}^{\,\ppm} $ at ${\bf X}_k$ by interpolation techniques; see \cite{Hallin2020DistributionAQ} for furter details.\medskip

These samples were generated from\vspace{-2mm}
\begin{enumerate}
\item[(i)] the centered bivariate normal distribution with covariance matrix
 $
\Sigma
=
\begin{pmatrix} 
	2 & 1\\ 
	1 & 1
\end{pmatrix}
$;
\item[(ii)]  
the bivariate distribution with independent exponential  marginals  with mean one;
\item[(iii)] the standard skew-$t$ distribution with~four degrees of freedom and slant \linebreak vector~$\alpha=(10,10)$--- see \cite{AzzCap2014}---and 
\item[(iv)] the banana-shaped  distribution considered in \cite{Hallin2020DistributionAQ}, viz. the mixture 
$$
\frac38 \NN(\mu_1, \Sigma_1) + \frac38\NN(\mu_2,\Sigma_2) + \frac14 \NN(\mu_3, \Sigma_3)
,
$$
of three  bivariate normal distributions
 where 
$$
\mu_1 = \begin{pmatrix}-3 \\ 0\end{pmatrix}
,\quad 
\mu_2 = \begin{pmatrix}3 \\ 0\end{pmatrix}
,\quad 
\mu_3 = \begin{pmatrix}0 \\ -\frac52\end{pmatrix}
,
$$
and 
$$
\Sigma_1= \begin{pmatrix} 5 & -4\\-4 & 5\end{pmatrix}
,\quad 
\Sigma_2 = \begin{pmatrix} 5 & 4\\4 & 5\end{pmatrix}
,\quad 
\Sigma_3 = \begin{pmatrix} 4 & 0\\0 & 1\end{pmatrix}
.
$$
\end{enumerate}

\begin{figure}\label{Fig1}
	\hspace*{-3mm}\includegraphics[width=1.04\textwidth]{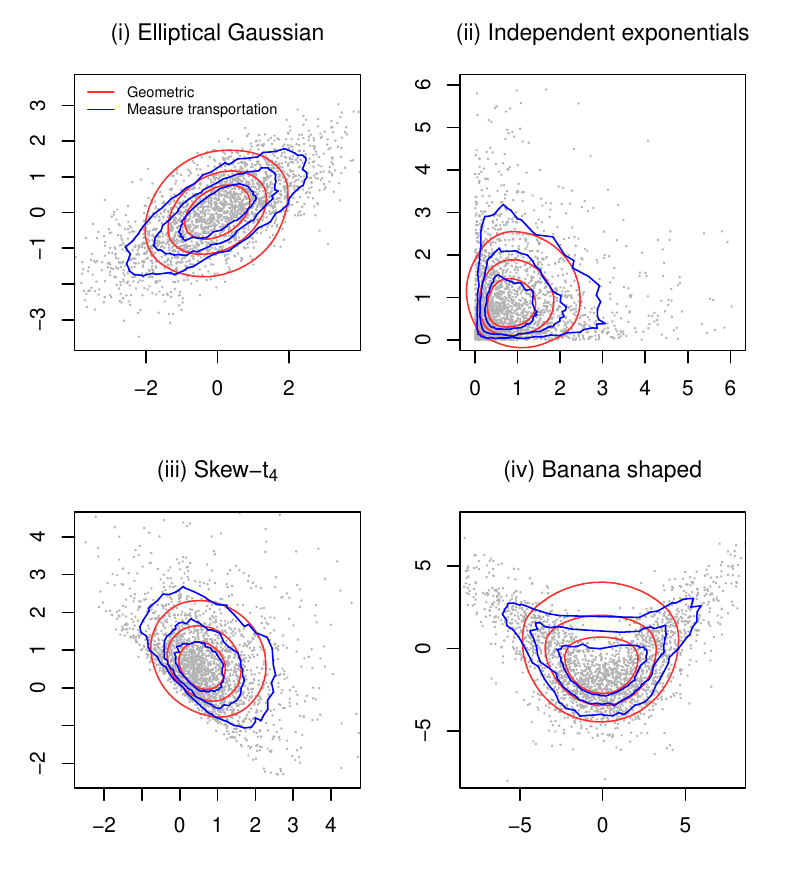}
	\caption{Geometric (red) and measure-transportation-based (blue) contours of levels $\tau=.25, .5, .75$ for samples of size $N=2400$ drawn from four different distributions. Geometric contours are relabeled according to their probability content, and measure-transportation contours are obtained through the interpolation scheme (see Hallin et al.~(2021)) of the optimal coupling between the observations and a $40\times60$ regular grid ${\mathfrak G}_N$ of the unit ball in $\R^2$; see Section 4.2 for details. Empirical quantiles, both geometric and center-outward, are connected by lines in the plots.
	}
\end{figure}\medskip

{Inspection of  Figure~1 reveals the following facts. 

\begin{enumerate}
\item[(a)]First, the contours obtained by measure transportation seem to better reflect the geometry of the underlying distributions. This is particularly obvious for the banana-shaped distribution (iv), where the concave shape of the mixture is well picked up by the center-outward contours but completely missed by the relabeled geometric ones. The same phenomenon also takes place, to a less spectacular extent, under  distributions (i)-(iii):  the relabeled geometric quantile region of order~.75, under the exponential distribution (ii), largely extends beyond the exponential support, something the center-outward contours do not; the same geometric  contours, for the skew-$t$ distribution~(iii), contrary to the center-outward ones,   hardly exhibit any skewness. 
\item[(b)]Second, while relabeled geometric quantiles (as a direct consequence of their re-labelization rather than because of  their L$_1$ characterization) do provide the desired control over the probability content of their regions, hence yield, for each~${\bf X}_i \in\X _N$, a sound concept $R^{(N)}_i\!\!\!\coloneqq \text{rank}(\Vert {\bf F}_{{\rm P}_N}^{\g}({\bf X}_i)\Vert )$ of center-outward multivariate rank, they do not provide any  clear notion of signs. Such signs, quite on the contrary, are well defined, in the measure-transportation-based approach,  for the sample points~${\bf X}_i \in~\!\X _n$, as the unit vectors ${\bf F}_{{\rm P}_N}^{\ppm}({\bf X}_i)/\Vert {\bf F}_{{\rm P}_N}^{\ppm}({\bf X}_i)\Vert$, $i=1,\ldots,N$, which are i.i.d. (uniform over the unit vectors associated with the $N$ gridpoints in ${\mathfrak G}_N$ and independent of the ranks $R^{(N)}_1 ,\ldots, R^{(N)}_N$). 
\item[(c)]On the other hand, geometric contours are more regular than their measure-transportation-based counterparts. This is a consequence of (iii) and (iv) in Sections 3.1 and~3.2, respectively: geometric contours are continuous as soon as ${\rm P}$ is non-atomic and smooth (at least $\CC^{d-1}$) while the oscillations, hence also the variability, of the interpolated measure-transportation-based  contours are further increased by the interpolation scheme. 
\end{enumerate}

While points (a) and (b) above are pleading for the measure-transportation-based  concepts, the decisive argument in favor of the latter follows from a theoretical result by  \cite{Girard17} on the ``intriguing behavior'' (a clear understatement) of high-order geometric quantile contours.

\subsection{Numerical comparison of extreme contours}

 \cite{Girard17} establish the following disturbing result for extreme geometric quantiles: if ${\rm P}$ is a probability measure on $\R^d$ with finite second-order moments, i.e. ${\rm E}_{\rm P}[\|{\bf X}\|^2]<\infty$ for ${\bf X}\sim {\rm P}$ with covariance matrix $\boldsymbol\Sigma$, then, denoting by~Tr$(\bf A)$ the trace of a square matrix $\bf A$,  
$$
\|\QPg(\tau {\bf u})\|^2 (1-\tau)
\to 
\frac{1}{2}\Big(
\Tr{\boldsymbol\Sigma}-{\bf u}^\prime {\boldsymbol\Sigma} {\bf u}
\Big)
$$
as $\tau\to 1$ for any ${\bf u}\in\S^{d-1}$. First, let us stress that the quantity $\Tr{\boldsymbol\Sigma}-{\bf u}^\prime {\boldsymbol\Sigma} {\bf u}$ is always non-negative, and is strictly positive as soon as $\boldsymbol\Sigma$ is non-degenerate. In addition to providing the rate at which extreme geometric quantiles are  ``escaping to infi\-nity,''  this result also provides valuable insights on the shape of the corres\-pon\-ding extreme contours. For instance, taking for~${\bf u}$  a unit eigenvector ${\bf u}_{\text{\rm max}}$ asso\-ciated~with~$\boldsymbol\Sigma$'s {\it largest} eigenvalue  yields a norm $\|\QPg(\tau {\bf u}_{\text{\rm max}})\|$ for the contour $\pmb{\CCC}_{\rm P}^{\g}(\tau)$ of order $\tau$ in  direction $\bf u$ which is   {\it minimal} for $\tau\to 1$;  the same  contour, in the direction~${\bf u}_{\text{\rm min}}$ associated with the smallest eigenvalue of $\boldsymbol\Sigma$,  asymptotically yields the {\it largest} norm~$\|\QPg(\tau {\bf u}_{\text{\rm min}})\|$. More generally, due to the fact that the shape factor~$\Tr{\boldsymbol\Sigma}-{\bf u}^\prime {\boldsymbol\Sigma} {\bf u}$ is a decreasing function of ${\bf u}^\prime {\boldsymbol\Sigma} {\bf u}$, the asymptotic behavior  of extreme quantiles---moving away fast from the center  along  the eigendirections of $\boldsymbol\Sigma$ associated with small eigenvalues, and  slowly along  the eigendirections of $\boldsymbol\Sigma$ associated with large eigenvalues---is exactly the opposite of what it is expected to be. This is highly counterintuitive and, actually,  totally unacceptable. 

Figure~2 provides a finite-sample empirical   illustration of this theoretical fact. We generated a sample of size $N=1,000$  from a centered bivariate  Gaussian distribution~$\NN({\boldsymbol 0},{\boldsymbol \Sigma})$ with covariance matrix 
 $
\boldsymbol \Sigma
=
\begin{pmatrix} 1/8 & 0\\0 & 3/4\end{pmatrix}
.
$ 
Then, for  $\tau\in\{.90, .95, .99\}$,  we computed 
the empirical relabeled geometric   and transportation-based quantile contours of order $\tau$, as described in Sections~\ref{41} and~\ref{42}, with $K=50$, $n_R=20$, and~$n_S=50$  for the latter (see Section \ref{42} for  notation).

\begin{figure}\label{Fig2}
	\begin{center}
	\includegraphics[width=0.6\textwidth]{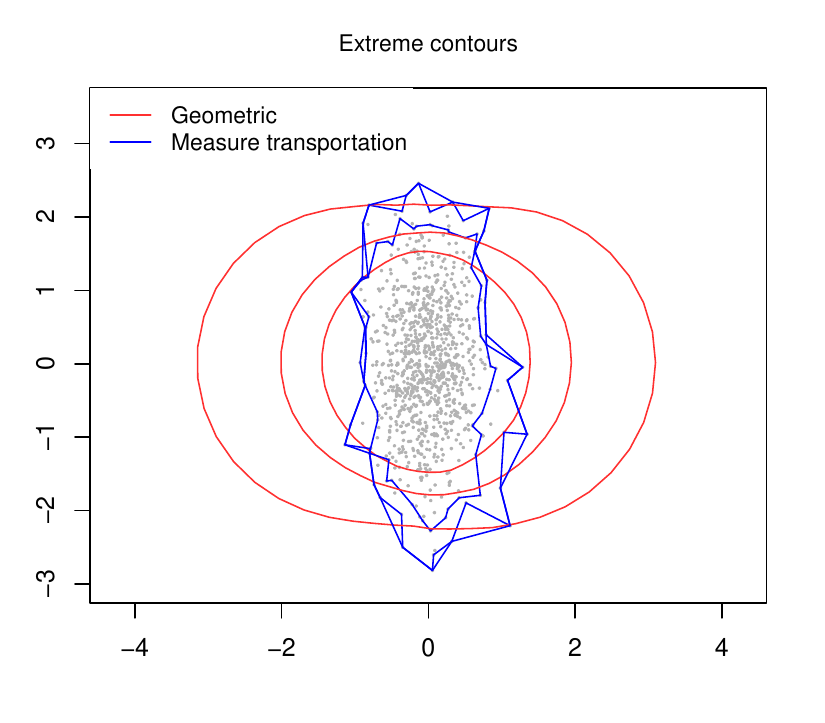}
	\end{center}
	\caption{Geometric (red) and measure-transportation-based (blue) contours of order  $\tau=.90, .95, .99$ for samples of size $N=1,000$ drawn from a (nonspherical)  bivariate Gaussian distribution on $\R^2$. Empirical quantiles, both geometric and center-outward, are connected by lines in the plots.
	}
\end{figure}

Figure 2 perfectly illustrates the pathological phenomenon described by  \cite{Girard17}: the geometric contours are much wider in the (horizontal) directions associated with $\boldsymbol\Sigma$'s   smaller eigenvalue, and narrower in the (vertical) directions associated with $\boldsymbol\Sigma$'s   larger eigenvalue. This pathological behaviour is increasingly severe as  $\tau$ approaches one:  for    $\tau = .90$, the contour is approximately circular. In sharp contrast, the center-outward contours, irrespective of $\tau$, are correctly accounting for  the shape of the distribution and the relative magnitudes of $\boldsymbol\Sigma$'s   eigenvalues.

\section{Conclusions}

Among all concepts of multivariate quantiles proposed in the literature, geometric quantiles (after due re-labelization) and measure-transportation-based center-outward quantiles are the most convincing ones: they both reduce to classical concepts in dimension one, and   both satisfy most properties expected from quantiles, quantile regions, and quantile contours. Empirical comparisons and a highly pathological property of extreme geometric quantiles established by  \cite{Girard17}, however, are tilting the balance   in favour of the measure-transportation-based  concept. 

\section*{Acknowledgments} 

Marc Hallin acknowledges the support of the Czech Science Foundation grant GA\v{C}R22036365. Dimitri Konen is supported by a research fellowship from the Centre for Research in Statistical Methodology (CRiSM) of the University of Warwick. We thank Eustasio del Barrio and Alberto Gonz{\'a}lez-Sanz for kindly providing their codes for the computation of center-outward quantile functions and contours.

\bibliographystyle{chicago}
\bibliography{Biblioreg}
\end{document}